\documentclass[11pt]{article}
\usepackage[a4paper,margin=27mm]{geometry}
\usepackage{amsmath,amssymb,amsthm}
\usepackage{hyperref}
\usepackage{authblk}
\hypersetup{colorlinks=true,linkcolor=blue,citecolor=blue,urlcolor=blue}
\newcommand{\PP}{\mathbb P}
\newcommand{\EE}{\mathbb E}
\newcommand{\ind}{\mathbf 1}
\newcommand{\GG}{\mathcal G}
\newtheorem{theorem}{Theorem}[section]

\newtheorem{lemma}{Lemma}[section]

\newtheorem{remark}{Remark}[section]
\numberwithin{equation}{section}
\newcommand{\msc}[1]{\par\noindent\textit{MSC 2020:} #1}

\title{Precise large deviations for the elephant random walk\\
with Weibull-type and lognormal-type random step sizes}
\author[1]{Shuyuan Liang\thanks{Corresponding author. E-mail: liangshuyuan0501@foxmail.com}}
\affil[1]{\footnotesize\itshape School of Mathematics and Statistics, Northeastern University at Qinhuangdao, Qinhuangdao, 066004, Hebei, China}
\date{}

\begin{document}
\maketitle

\begin{abstract}
We prove precise large deviations for the symmetric elephant random walk with random step sizes. When the distribution of \(|Z_1-1|\) is subexponential, we first show that, for each fixed \(p\in[0,1]\), there exists a positive number sequence \(b_n\to\infty\) such that
\[
\PP(S_n>x)\sim \PP(S_n<-x)\sim \frac n2\overline F(x)
\]
uniformly for \(x\ge b_n\), where \(\overline F(x)=\PP(|Z_1-1|>x)\). This existence result requires no moment assumptions and does not specify the growth rate of \(b_n\). For \(p\in[0,3/4]\), we further give explicit sufficient conditions on \(b_n\) under Weibull-type and lognormal-type tail conditions with a regularly varying factor \(x^\gamma L(x)\), where \(\gamma\in\mathbf R\) and \(L\) is a positive measurable slowly varying function.  Analogous results hold with \(\PP(|Z_1|>x)\) and non-strict tails.
\end{abstract}

\msc{60F10; 60G50}
\section{Introduction}\quad
The elephant random walk (ERW) is a one-dimensional random walk model with entire history memory defined on $\mathbf{Z}$. In 2004, Schütz and Trimper \cite{schutz2004elephants} first introduced this model arising from mathematical physics to study how the long-term memory has an impact on non-Markovian random walks. The definition of the ERW is as follows.  
Let \((X_n)_{n\ge1}\) and $(T_n)_{n\geq0}$ be respectively the step and the position of the ERW at time $n$ and the following equalities show the connection between $(X_n)_{n\geq1}$ and $(T_n)_{n\geq0}$:
\begin{equation*}
\begin{cases}
T_0=0,\\
T_n=T_{n-1}+X_n,\quad&n\geq1,
\end{cases}
\end{equation*}
which implies that $T_n=\sum_{k=1}^{n}X_k$, $n\geq1$.
For all $n\geq1$, $X_n=s\in\{1,-1\}$.  Set $q\in(0,1)$, memory parameter $p\in[0,1]$ and \(\GG_{n}=\sigma(X_1,\ldots,X_{n})\), thus the distribution of $X_n$ is
\begin{equation*}
\begin{cases}
\mathbb{P}(X_1=s)=sq+\mathbf1_{\{s=-1\}}\\
\mathbb{P}(X_n=s|\GG_{n-1})=\dfrac12+\dfrac{s(2p-1)}{2(n-1)}T_{n-1},\quad&n\geq2.\\
\end{cases}
\end{equation*}\par
For the ERW model, many properties have been studied, especially in asymptotic behaviors. Coletti \textit{et al.} \cite{CGSCLT} have obtained the central limit theorem (CLT) for the ERW. Using martingale method, Bercu \cite{BMERW} also derived the CLT results.  In \cite{B21}, Bertoin studied how the memory parameter $p$ impacts the number of zeros when $0\leq p<3/4$. In \cite{FanHuMa2021}, Fan \textit{et al.} have established the Cram\'er-type moderate deviations for the ERW.
\par
In this paper, we consider a more generalized ERW model with random step sizes, which was first introduced by Fan and Shao \cite{FSCMD}. 
Let \((Z_n)_{n\ge1}\) be an i.i.d. real-valued random variable sequence, independent of the entire step
sequence $(X_n)_{n\geq1}$. Let $S_n$ be the position of the ERW with random step sizes at time $n$, that is,
\[
S_n=\sum_{k=1}^nX_kZ_k.
\]
When \(\EE Z_1=\nu\neq0\) and \(\operatorname{Var}(Z_1)<\infty\), we have
\[
S_n=\nu T_n+\sum_{k=1}^nX_k(Z_k-\nu).
\]
Replacing $Z_n$ by $Z_n/\nu$ and $S_n$ by $S_n/\nu$ allows us to assume $\nu=1$. We write $\sigma^2=\operatorname{Var}(Z_1)$ for the variance after normalization. When $\sigma=0$, the ERW with random step sizes degenerates to the general ERW model.\par
Recently, there are many works about how the distribution of $Z_1$ affects the properties of the ERW with random step sizes. In \cite{FSCMD}, Fan and Shao have established the Cram\'er-type moderate deviations when $Z_1$ satisfies Bernstein's condition. When Cram\'er's condition is not satisfied, Fan \textit{et al.} \cite{FanLiangMa} obtained the Cram\'er-type moderate deviations and precise large deviations for the ERW with random step sizes. With local regularly varying tails, Liang \cite{Liang2026} derived local precise large deviations.\par
Our paper studies precise large deviations for the symmetric ERW ($q=1/2$) with random step sizes. When the distribution of \(|Z_1-1|\) is subexponential, Theorem \ref{thm:subexponential} first establishes the existence of a positive number sequence \(b_n\to\infty\) for which both tails of \(S_n\) are uniformly asymptotic to \(n\PP(|Z_1-1|>x)/2\) for \(x\ge b_n\), for each fixed \(p\in[0,1]\). This result follows from the fixed-sum property of subexponential distributions and the sign symmetry of the ERW, but does not give an explicit growth rate of \(b_n\). We then consider Weibull-type and lognormal-type tail conditions with a regularly varying factor. Theorems \ref{thm:main} and \ref{thm:lognormal} give uniform asymptotic estimates for both tails of \(S_n\) when \(x\ge b_n\), with explicit sufficient conditions on \(b_n\) for \(p\in[0,3/4]\). For independent sums, related results for these tail classes are given in Theorems 8.2 (ii) and 8.3 of Denisov \textit{et al.} \cite{DDS}. Here we give preliminary estimates that are uniform over deterministic sign vectors, and then use conditioning and an Azuma--Hoeffding-type inequality to treat the ERW. Section 2 states the main results, Section 3 gives preliminary lemmas, and Section 4 proves the main results, with Subsections 4.1, 4.2, and 4.3 devoted to the general subexponential, Weibull-type, and lognormal-type results, respectively.\par
Throughout the paper, we denote by $C$ the universal positive constant and $C_\tau$ the positive constant only depending on $\tau$. Besides, their value may vary in different expressions. For two positive sequences $(p_n)_{n\geq1}$ and $(q_n)_{n\geq1}$, the notations $p_n\gg q_n$ and $p_n\ll q_n$ mean that $p_n/q_n\to\infty$ and $p_n/q_n\to0$ respectively. For two positive functions $\mathcal{F}(x)$ and $\mathcal{G}(x)$, use $\mathcal{F}(x)\sim\mathcal{G}(x)$ to denote $\mathcal{F}(x)/\mathcal{G}(x)\to1$ as $x\to\infty$.
\section{Main results}
\quad Let \(F\) denote the
distribution of \(|Z_1-1|\), with
\(\overline F(t)=\PP(|Z_1-1|>t)\). We first consider the case when
\(F\) is subexponential. More precisely, we write \(F\in\mathcal S\)
if \(\overline F(t)>0\) for all \(t>0\) and
\[
\PP\!\left(\sum_{k=1}^m|Z_k-1|>t\right)
\sim m\overline F(t),
\]
where $m\geq2$
and \(F\) is long-tailed, that is,
\(\overline F(t+a)/\overline F(t)\to1\) for every fixed
\(a\in\mathbf R\). These are standard properties of subexponential
distributions; see Nagaev A V \cite{NAVILTLDP}, Nagaev S V \cite{NSVLDP} and Denisov \textit{et al.} \cite{DDS}.

\begin{theorem}\label{thm:subexponential}
Assume that \(F\in\mathcal S\), with no moment assumptions on \(Z_1\). Then, for \(q=1/2\) and each fixed
\(p\in[0,1]\), there exists a positive number sequence
\((b_n)_{n\ge1}\) with \(b_n\to\infty\) such that
\begin{equation}\label{subexp:main-right}
\lim_{n\to\infty}\sup_{x\ge b_n}
\left|\frac{\PP(S_n>x)}{n\overline F(x)}-\frac12\right|=0
\end{equation}
and
\begin{equation}\label{subexp:main-left}
\lim_{n\to\infty}\sup_{x\ge b_n}
\left|\frac{\PP(S_n<-x)}{n\overline F(x)}-\frac12\right|=0.
\end{equation}
The same conclusions hold with \(\PP(|Z_1|>x)\) in place of
\(\overline F(x)\), and with non-strict tail inequalities.
\end{theorem}

\par

Theorem \ref{thm:subexponential} gives the existence of a sufficiently
large threshold under the general subexponential condition. However,
it does not specify the growth rate of \(b_n\), or assert that
\(b_n\) can be taken of linear order. To obtain explicit sufficient
conditions on \(b_n\), we next consider two particular subexponential
tail classes, namely Weibull-type and lognormal-type tails with a
regularly varying factor. For these two classes, assume that \(\EE Z_1=1\)
and write \(\sigma^2=\operatorname{Var}(Z_1)\); the tail conditions below imply that \(\sigma^2<\infty\).
Let \(\gamma\in\mathbf R\) be fixed and \(c_1,c_{-1}\ge0\) with \(c_1+c_{-1}>0\). Set \(L\) as
 a positive measurable function slowly varying at infinity, that is,
\[
\frac{L(rt)}{L(t)}\to1\qquad\textrm{as }t\to\infty
\quad\textrm{for every }r>0,
\]
Thus \(t^\gamma L(t)\) is regularly varying with index \(\gamma\).
The case \(L\equiv1\) gives a power factor, and \(\gamma=0\) gives a
slowly varying factor.
In each result below, \(H\) denotes the corresponding tail comparison function;
\(H\) itself is not required to be a distribution tail.

We first consider Weibull-type tails. Let \(0<\alpha<1\), \(\beta>0\), and set
\begin{equation}\label{H:weibull}
H(t)=t^\gamma L(t)\exp\{-\beta t^\alpha\},\qquad t>0.
\end{equation}
Assume that the following tail condition holds:
\begin{equation}\label{WeiTail}
\frac{\PP(s(Z_1-1)>t)}{H(t)}\to c_s
\qquad\textrm{as }t\to\infty,\qquad s\in\{-1,1\}.
\end{equation}
Then we have \(\overline F(t)\sim(c_1+c_{-1})H(t)\). When $Z_1$ satisfies the condition \eqref{WeiTail}, we have the following theorem.

\begin{theorem}\label{thm:main}
Under the preceding condition \eqref{WeiTail}, let 
\((b_n)_{n\geq1}\) 
be a positive number sequence satisfying 
\begin{equation}\label{proof:boundary}
\begin{cases}
b_n\gg n^{1/(2-2\alpha)},\quad&p\in[0, 3/4)\\
b_n\gg (n\log n)^{1/(2-2\alpha)},\quad&p=3/4.
\end{cases}
\end{equation}
Then, for $q=1/2$ and each fixed $p\in[0,3/4]$, it holds that
\begin{equation}\label{eq:main-right}
\lim_{n\to\infty}\sup_{x\ge b_n}
\left|\frac{\PP(S_n>x)}{n\overline F(x)}-\frac12\right|=0
\end{equation}
and
\begin{equation}\label{eq:main-left}
\lim_{n\to\infty}\sup_{x\ge b_n}
\left|\frac{\PP(S_n<-x)}{n\overline F(x)}-\frac12\right|=0.
\end{equation}
The same conclusions hold with \(\PP(|Z_1|>x)\) in place of
\(\overline F(x)\), and with non-strict tail inequalities.
\end{theorem}\par

Next, we consider lognormal-type tails. Write
\[
\overline\Phi(z)=\frac{1}{\sqrt{2\pi}}\int_z^\infty e^{-u^2/2}\,du,
\qquad z\in\mathbf R.
\]
For \(\beta>0\) and \(\mu\in\mathbf R\), set
\begin{equation}\label{H:lognormal}
H(t)=t^\gamma L(t)
\overline\Phi\!\left(\sqrt{2\beta}(\log t-\mu)\right),\qquad t>0.
\end{equation}
The lognormal factor has logarithmic mean \(\mu\) and logarithmic
variance \(1/(2\beta)\). Assume that
\begin{equation}\label{LogTail}
\frac{\PP(s(Z_1-1)>t)}{H(t)}\to c_s
\qquad\textrm{as }t\to\infty,\qquad s\in\{-1,1\},
\end{equation}
where \(H\) is now given by \eqref{H:lognormal}. Again, it holds that
\(\overline F(t)\sim(c_1+c_{-1})H(t)\).

\begin{theorem}\label{thm:lognormal}
Under condition \eqref{LogTail}, let \((b_n)_{n\ge1}\) be a positive
number sequence satisfying
\begin{equation}\label{log:boundary}
\begin{cases}
b_n\gg\sqrt n\,\log n,\quad&0\le p<3/4,\\
b_n\gg\sqrt{n\log n}\,\log n,\quad&p=3/4.
\end{cases}
\end{equation}
Then, for \(q=1/2\) and each fixed \(p\in[0,3/4]\), the conclusions
\eqref{eq:main-right} and \eqref{eq:main-left} hold. The same conclusions hold with
\(\PP(|Z_1|>x)\) in place of \(\overline F(x)\), and with non-strict
tail inequalities.
\end{theorem}

\begin{remark}\label{rem:lognormal-form}
By the normal-tail asymptotic,
\[
\overline\Phi\!\left(\sqrt{2\beta}(\log x-\mu)\right)
\sim\frac{\exp\{-\beta(\log x-\mu)^2\}}
{2\sqrt{\pi\beta}(\log x-\mu)}.
\]
In particular, Theorem \ref{thm:lognormal} also covers tail comparison functions
of the form
\[
H(x)=x^\gamma L(x)\exp\{-\beta(\log x)^2\}.
\]
Indeed, take \(\mu=0\) and replace the slowly varying factor in \eqref{H:lognormal}
by \(2\sqrt{\pi\beta}(\log x)L(x)\) for all sufficiently large \(x\).
The two tail comparison functions are then asymptotically equivalent, so the
conclusions and the sufficient boundaries remain the same.
The tail limits in \eqref{WeiTail} and \eqref{LogTail} allow either
\(c_1\) or \(c_{-1}\) to vanish; in particular, positive random step
sizes are included.
\end{remark}

\section{Preliminary lemmas}
\quad
In this section, \(H\) is given by \eqref{H:weibull} under
\eqref{WeiTail}, and by \eqref{H:lognormal} under \eqref{LogTail}.
We first recall some properties of slowly varying functions. By the
uniform convergence theorem and Potter's bound (see Theorems 1.2.1 and
1.5.6 in Bingham \textit{et al.} \cite{BGT}), we have
\begin{equation}\label{slow:uniform}
\sup_{r\in[a,b]}\left|\frac{L(rx)}{L(x)}-1\right|\to0
\quad\textrm{as }x\to\infty,\qquad 0<a<b<\infty,
\end{equation}
and, for every \(\eta>0\), there exists a positive constant $C_\eta$ such that 
\begin{equation}\label{slow:growth}
C_\eta^{-1}x^{-\eta}\le L(x)\le C_\eta x^\eta
\quad\textrm{for all sufficiently large }x.
\end{equation}
Consequently,
\begin{equation}\label{rv:uniform}
\sup_{r\in[a,b]}
\left|\frac{(rx)^\gamma L(rx)}{x^\gamma L(x)}-r^\gamma\right|\to0,
\qquad 0<a<b<\infty.
\end{equation}
Moreover, for every \(\varepsilon>0\), there exists \(x_0>0\) such that
\begin{equation}\label{rv:potter}
\left|\log\frac{y^\gamma L(y)}{x^\gamma L(x)}\right|
\le\log(1+\varepsilon)+(|\gamma|+1)\left|\log\frac yx\right|,
\qquad x,y\ge x_0.
\end{equation}
Indeed, apply Potter's bound to \(L\) with exponent \(1\) and constant
\(1+\varepsilon\), interchange \(x\) and \(y\), and include the
factor \((y/x)^\gamma\).
These estimates do not require differentiability or monotonicity of \(L\).

For \(s\in\{-1,1\}\), define
\(F_s(t)=\PP(s(Z_1-1)>t)\). Define the nonincreasing error function
\begin{equation}\label{rho}
\rho(u)=\max_{s\in\{-1,1\}}\sup_{t\ge u}
\left|\frac{F_s(t)}{H(t)}-c_s\right|.
\end{equation}
Then we have \(\rho(u)\to0\) as \(u\to\infty\), and
\(F_s(t)\le C H(t)\) for all sufficiently large \(t\), uniformly in \(s\).

\subsection{Weibull-type tails}
\quad Under condition \eqref{WeiTail}, for all \(0<\delta<\beta\), it holds that
\begin{equation}
\EE\!\left[(Z_1-1)^2
\exp\{(\beta-\delta)|Z_1-1|^\alpha\}\right]<\infty.
\label{eq:weighted-moment}
\end{equation}
Indeed, by \eqref{slow:growth}, for every fixed \(\eta>0\) and all
sufficiently large \(t\),
\(\PP(|Z_1-1|>t)\le Ct^{\gamma+\eta}\exp\{-\beta t^\alpha\}\).
Integration by parts reduces the tail part of the expectation in
\eqref{eq:weighted-moment} to an integral bounded by %a constant times
\(\int_1^\infty(t^{1+\gamma+\eta}+t^{1+\alpha+\gamma+\eta})e^{-\delta t^\alpha}\,dt\),
which is finite.
Therefore, we can obtain the following upper bound of the truncated exponential moment.
\begin{lemma}\label{lem:mgf}
Fix \(0<\delta<\beta\). For \(s\in\{-1,1\}\), \(y>0\),
\(\lambda\ge0\), and \(\lambda y^{1-\alpha}\le\beta-\delta\), we have the following inequality
\begin{equation*}
\EE\!\left[\exp\{\lambda s(Z_1-1)\}\ind_{\{s(Z_1-1)\le y\}}\right]
\le1+\frac{\lambda^2}{2}
\EE\!\left[(Z_1-1)^2\exp\{(\beta-\delta)|Z_1-1|^\alpha\}\right].
\end{equation*}
\end{lemma}

\begin{proof}
For every \(z\in\mathbf R\), we have
\begin{equation}\label{expz}
e^z\le1+z+\frac{z^2}{2}e^{z_+},
\end{equation}
where $z_+=\max(z,0).$
On \(0\le s(Z_1-1)\le y\) and by \(\lambda y^{1-\alpha}\le\beta-\delta\), it holds that
\begin{equation}\label{alpha0}
\lambda s(Z_1-1)\le(\beta-\delta)[s(Z_1-1)]^\alpha.
\end{equation}
Since \(\EE[s(Z_1-1)]=0\) and \(y>0\), we have
\begin{equation}\label{Es}
\EE\!\left[s(Z_1-1)\ind_{\{s(Z_1-1)\le y\}}\right]
=-\EE\!\left[s(Z_1-1)\ind_{\{s(Z_1-1)>y\}}\right]\le0.
\end{equation}
Consequently, using \eqref{expz}-\eqref{Es}, the following inequalities hold
\begin{align*}
&\quad\EE\!\left[\exp\{\lambda s(Z_1-1)\}\ind_{\{s(Z_1-1)\le y\}}\right]\\
&\le \PP(s(Z_1-1)\le y)
+\lambda\EE\!\left[s(Z_1-1)\ind_{\{s(Z_1-1)\le y\}}\right]\\
&\quad+\frac{\lambda^2}{2}\EE\!\left[
  [s(Z_1-1)]^2\exp\{\lambda[s(Z_1-1)]_+\}\ind_{\{s(Z_1-1)\le y\}}\right]\\
&\le1+\frac{\lambda^2}{2}\EE\!\left[
  [s(Z_1-1)]^2\exp\{(\beta-\delta)|s(Z_1-1)|^\alpha\}\right]\\
&=1+\frac{\lambda^2}{2}
\EE\!\left[(Z_1-1)^2\exp\{(\beta-\delta)|Z_1-1|^\alpha\}\right].
\end{align*}
Hence, we complete the proof of Lemma \ref{lem:mgf}.
\end{proof}\par
The following lemma states the precise large deviation for $\PP\!\left(\sum_{k=1}^ns_k(Z_k-1)>x\right)$, where $s_k\in\{1,-1\}$.
\begin{lemma}\label{lem:uniform}
Assume that $Z_1$ satisfies the condition \eqref{WeiTail}, and let
\((b_n)_{n\ge1}\) satisfy \(b_n\gg n^{1/(2-2\alpha)}\). For every
deterministic sign vector
\(\boldsymbol s=(s_1,\ldots,s_n)\in\{-1,1\}^n\), we have
\begin{equation*}
\lim_{n\to\infty}
\sup_{\boldsymbol s\in\{-1,1\}^n}\sup_{x\ge b_n}
\left|\frac{\PP\!\left(\sum_{k=1}^ns_k(Z_k-1)>x\right)}{nH(x)}
-\frac1n\sum_{k=1}^nc_{s_k}\right|=0.
\end{equation*}
\end{lemma}

\begin{proof}
All estimates in this proof are uniform in \(\boldsymbol s\) and \(x\ge b_n\).
The variables \((s_k(Z_k-1))_{1\leq k\leq n}\) are independent, centered, and have
variance \(\sigma^2\).
Condition $b_n\gg n^{1/(2-2\alpha)}$ implies that
\begin{equation}
n b_n^{-2(1-\alpha)}\to0,\qquad
b_n^\alpha/\log n\to\infty.
\label{eq:threshold-consequences}
\end{equation}
Choose constants $\theta$ and $\kappa$ such that
\begin{equation}
2^{-1/\alpha}<\theta<\frac12,\qquad
\beta<\kappa<\beta\theta^{\alpha-1}.
\label{eq:constants}
\end{equation}
These choices are possible because \(0<\alpha<1\).
Fix \(x\ge b_n\).
There is a fixed \(0<\delta<\beta\) such that
\[
\kappa\theta^{1-\alpha}\le\beta-2\delta,\qquad
\beta\alpha\le\beta-2\delta.
\]
For all sufficiently large \(n\), uniformly for \(x\ge b_n\), we have
\(\lambda (\theta x)^{1-\alpha}\le\beta-\delta\) for either
\(\lambda=\kappa x^{\alpha-1}\) or
\(\lambda=\beta\alpha (\theta x)^{\alpha-1}+|\gamma|/(\theta x)\).
Indeed, for the second choice the left-hand side equals
\(\beta\alpha+|\gamma|(\theta x)^{-\alpha}\).
Thus Lemma \ref{lem:mgf} applies to both choices. For either choice, it holds that
\begin{equation}
\sup_{x\ge b_n}n\lambda^2\to0,\qquad
\sup_{x\ge b_n}\sqrt n\,\lambda\to0.
\label{eq:lambda-small}
\end{equation}\par
According to the number of terms satisfying $s_k(Z_k-1)>\theta x$, $\PP\!\left(\sum_{k=1}^ns_k(Z_k-1)>x\right)$ can be decomposed as the following three parts.
\begin{align}\label{Z_kPLDdec}
&\quad\PP\!\left(\sum_{k=1}^ns_k(Z_k-1)>x\right)\notag\\
&=\PP\left(\sum_{k=1}^ns_k(Z_k-1)>x,\ \max_{1\leq k\leq n}s_k(Z_k-1)\le \theta x\right)\notag\\
&\quad+\sum_{i=1}^{n}\PP\left(\sum_{k=1}^ns_k(Z_k-1)>x,\ s_i(Z_i-1)>\theta x,\ \max_{j\neq i}s_j(Z_j-1)\leq\theta x\right)\notag\\
&\quad+\PP\left(\sum_{k=1}^ns_k(Z_k-1)>x,\ \#\{k:s_k(Z_k-1)>\theta x\}\ge2\right)\notag\\
&=:P_{1n}+P_{2n}+P_{3n}.
\end{align}
Then we give the estimates of the upper bound for $P_{1n}$, $P_{2n}$ and $P_{3n}$ respectively.\par\noindent
\textbf{(I) Control of $P_{1n}.$} Take \(\lambda=\kappa x^{\alpha-1}\).
The exponential Markov inequality and Lemma \ref{lem:mgf} yield
\begin{align}\label{P1nup}
P_{1n}
&\le \exp\{-\lambda x\}
\prod_{k=1}^n\EE\!\left[\exp\{\lambda s_k(Z_k-1)\}\ind_{\{s_k(Z_k-1)\le\theta x\}}\right]\notag\\
&\le\exp\{-\kappa x^\alpha+n\ln(1+C\lambda^2)\}
\le\exp\{-\kappa x^\alpha+Cn\lambda^2\}
=o(nH(x)).
\end{align}
The last assertion follows from \(\kappa>\beta\),
\eqref{slow:growth}, and \eqref{eq:lambda-small}, since
\[
\frac{P_{1n}}{nH(x)}
\le\frac{1}{nx^\gamma L(x)}
\exp\{-(\kappa-\beta)x^\alpha+Cn\lambda^2\}=o(1).
\]
\par\noindent
\textbf{(II) Control of $P_{2n}$.}
For each \(i\), define
\[
R_i=\sum_{j\ne i}s_j(Z_j-1),\qquad
E_i=\left\{\max_{j\ne i}s_j(Z_j-1)\le \theta x\right\}.
\]
Independence of $(Z_n)_{n\geq1}$ gives the exact identity
\begin{equation}\label{eq:one-jump-identity}
\PP\left(\sum_{k=1}^{n}s_k(Z_k-1)>x,\ s_i(Z_i-1)>\theta x,\ E_i\right)
=\EE\!\left[
F_{s_i}\big(\max\{\theta x,x-R_i\}\big)\ind_{E_i}\right].
\end{equation}
For all \(r\in\mathbf R\), the mean value theorem implies that
\begin{equation}\label{deltaxalpha}
x^\alpha-\max\{\theta x,x-r\}^{\alpha}\leq\frac{\alpha\left(x-\max\{\theta x,x-r\}\right)}{\max\{\theta x,x-r\}^{1-\alpha}}
\le\alpha (\theta x)^{\alpha-1}r_+.
\end{equation}
By \eqref{rho}, we have
\begin{equation}\label{Fsiup}
F_{s_i}\left(\max\{\theta x,x-R_i\}\right)\leq\left(c_{s_i}+\rho(\theta x)\right)H\left(\max\{\theta x, x-R_i\}\right).
\end{equation}
On \(\{R_i\ge0\}\), the argument \(\max\{\theta x,x-R_i\}\)
lies in \([\theta x,x]\), and we can deduce that
\[
\left|\log\frac{\max\{\theta x,x-R_i\}}{x}\right|
\le\frac{(R_i)_+}{\theta x}.
\]
Thus, by \eqref{slow:uniform} and
\eqref{deltaxalpha}, it holds that
\begin{align}\label{Hup}
H\left(\max\{\theta x, x-R_i\}\right)
&=H(x)\left(\frac{\max\{\theta x,x-R_i\}}{x}\right)^\gamma
\frac{L(\max\{\theta x,x-R_i\})}{L(x)}\notag\\
&\quad\times\exp\left\{\beta\left(x^\alpha-\max\{\theta x, x-R_i\}^\alpha\right)\right\}\notag\\
&\leq(1+o(1))H(x)
\exp\left\{\left(\beta\alpha(\theta x)^{\alpha-1}
+\frac{|\gamma|}{\theta x}\right)(R_i)_+\right\},
\end{align}
where the error is deterministic and uniform in \(i\), \(x\ge b_n\),
and \(R_i\ge0\).
On \(\{R_i<0\}\), monotonicity of \(F_{s_i}\) gives instead
\[
F_{s_i}(x-R_i)\le F_{s_i}(x)
\le(c_{s_i}+\rho(\theta x))H(x).
\]
Taking \(\lambda=\beta\alpha (\theta x)^{\alpha-1}+|\gamma|/(\theta x)\) and combining
the two cases, we obtain the following inequality:
\begin{equation}\label{FsiEup}
\EE\!\left[
F_{s_i}\big(\max\{\theta x,x-R_i\}\big)\ind_{E_i}\right]
\le(c_{s_i}+\rho(\theta x))(1+o(1))H(x)
\EE\!\left[\exp\{\lambda(R_i)_+\}\ind_{E_i}\right].
\end{equation}
For all \(r\in\mathbf R\), write \(r_-=\max\{-r,0\}\). Using the inequality
\[
e^{\lambda r_+}\le e^{\lambda r}+\lambda r_-,
\qquad \lambda\ge0
\]
and independence of $(Z_n)_{n\geq1}$, we can imply that
\begin{align}\label{Eexpup}
\EE\!\left[\exp\{\lambda(R_i)_+\}\ind_{E_i}\right]
&\le\EE\!\left[\exp\{\lambda R_i\}\ind_{E_i}\right]
 +\lambda\EE(R_i)_-\notag\\
&\le \exp\{Cn\lambda^2\}
 +\lambda\sqrt{\EE R_i^2}\notag\\
&=\exp\{Cn\lambda^2\}
 +\lambda\sigma\sqrt{n-1}
=1+o(1).
\end{align}
By \eqref{eq:one-jump-identity}-\eqref{Eexpup}, we get the following upper bound:
\begin{equation*}
\PP\left(\sum_{k=1}^{n}s_k(Z_k-1)>x,\ s_i(Z_i-1)>\theta x,\ E_i\right)
\leq\left(c_{s_i}+\rho(\theta x)\right)H(x)(1+o(1)), 
\end{equation*}
which implies that
\begin{equation}\label{P2nup1}
P_{2n}\leq\left(\sum_{k=1}^{n}c_{s_k}+n\rho(\theta x)\right)H(x)(1+o(1)).
\end{equation}
As a fact of
\begin{equation*}
0\leq\sum_{k=1}^{n}c_{s_k}\leq\max\{c_1, c_{-1}\}n,
\end{equation*}
it holds that
\begin{equation}\label{sumcOn}
\sum_{k=1}^{n}c_{s_k}=O(n).
\end{equation}
For all $x\geq b_n$, we have
\begin{equation}\label{nrhoo1}
n\rho(\theta x)=o(n).
\end{equation}
Using \eqref{P2nup1}-\eqref{nrhoo1}, we get the upper bound of $P_{2n}$:
\begin{equation}\label{P2nup}
P_{2n}\leq H(x)\left(\sum_{k=1}^{n}c_{s_k}+o(n)\right).
\end{equation}
\par\noindent
\textbf{(III) Control of $P_{3n}$.}
Since $(Z_n)_{n\geq1}$ is an i.i.d random variable sequence, we can obtain that
\begin{align*}
P_{3n}&\leq\PP\big(\#\{k:s_k(Z_k-1)>\theta x\}\ge2\big)\notag\\
&\le\sum_{1\leq i<j\leq n}F_{s_i}(\theta x)F_{s_j}(\theta x)\notag\\
&\le Cn^2(\theta x)^{2\gamma}L(\theta x)^2\exp\{-2\beta\theta^\alpha x^\alpha\}.
\end{align*}
After division by \(nH(x)\), we have the following upper bound
\[
\frac{P_{3n}}{nH(x)}
\le Cn x^\gamma\frac{L(\theta x)^2}{L(x)}
\exp\{-\beta(2\theta^\alpha-1)x^\alpha\}=o(1).
\]
Here the last equality follows from \(2\theta^\alpha>1\),
\eqref{slow:uniform}, \eqref{slow:growth}, and
\eqref{eq:threshold-consequences}. This implies that
\begin{equation}\label{P3nup}
P_{3n}= o(nH(x)).
\end{equation}\par
Combining \eqref{P1nup}, \eqref{P2nup} and \eqref{P3nup}, it holds that
\begin{equation}\label{PLDup}
\PP\!\left(\sum_{k=1}^ns_k(Z_k-1)>x\right)
\leq H(x)\left(\sum_{k=1}^{n}c_{s_k}+o(n)\right).
\end{equation}\par
Next, we will give an estimate for the lower bound. For all $n\geq1$, let $d_n=n^{1/4}b_n^{(1-\alpha)/2}$. Therefore, it holds that 
\[
\sqrt n\ll d_n\ll b_n^{1-\alpha}.
\]
The events
\[
\{s_i(Z_i-1)>x+d_n,\ R_i\ge-d_n,\ E_i\},\qquad 1\leq i\leq n
\]
are disjoint and are all contained in \(\{\sum_{j=1}^ns_j(Z_j-1)>x\}\).
By Chebyshev's inequality and addition rule of probability, we can deduce that
\begin{align}\label{PR_ilow}
\PP(R_i\ge-d_n,E_i)
&=\PP(E_i)-\PP(R_i<-d_n, E_i)\notag\\
&\geq1-\PP\left(\bigcup_{j\neq i}\{s_j(Z_j-1)>\theta x\}\right)-\PP(R_i<-d_n)\notag\\
&\ge1-\sum_{j\ne i}F_{s_j}(\theta x)-\frac{(n-1)\sigma^2}{d_n^2}\notag\\
&\ge1-CnH(\theta x)-\frac{n\sigma^2}{d_n^2}=1-o(1).
\end{align}
Moreover, it holds that
\[
0\le\beta\big((x+d_n)^\alpha-x^\alpha\big)
\le\beta\alpha d_n x^{\alpha-1}
\le\beta\alpha d_n b_n^{\alpha-1}\to0.
\]
Also, \(d_n/x\le d_n/b_n\to0\), so \eqref{slow:uniform} yields
\(L(x+d_n)/L(x)\to1\) uniformly for \(x\ge b_n\).
Moreover, \((1+d_n/x)^\gamma\to1\) uniformly on this range. Consequently,
\[
\frac{H(x+d_n)}{H(x)}
=\left(1+\frac{d_n}{x}\right)^\gamma\frac{L(x+d_n)}{L(x)}
\exp\{-\beta((x+d_n)^\alpha-x^\alpha)\}\to1
\]
uniformly for \(x\ge b_n\), and
\begin{equation}\label{lontail}
\sum_{i=1}^nF_{s_i}(x+d_n)
=H(x)\left(\sum_{i=1}^nc_{s_i}+o(n)\right)
\end{equation}
uniformly in the signs as well. By \eqref{PR_ilow} and \eqref{lontail}, we can obtain the following lower bound:
\begin{align}\label{PLDlow}
\PP\left(\sum_{k=1}^ns_k(Z_k-1)>x\right)
&\ge\sum_{i=1}^n\PP(s_i(Z_i-1)>x+d_n,\ R_i\ge-d_n,\ E_i)\notag\\
&=\sum_{i=1}^nF_{s_i}(x+d_n)\PP(R_i\ge-d_n,E_i)\notag\\
&\ge H(x)\left(\sum_{k=1}^{n}c_{s_k}-o(n)\right).
\end{align}\par
Combining \eqref{PLDup} and \eqref{PLDlow}, we complete the proof of Lemma \ref{lem:uniform}.
\end{proof}
\subsection{Lognormal-type tails}
\quad We now assume \eqref{LogTail}, with \(H\) given by
\eqref{H:lognormal}. Integration by parts gives the normal-tail bounds
\begin{equation}\label{log:mills}
\frac{z}{1+z^2}\frac{e^{-z^2/2}}{\sqrt{2\pi}}
\le\overline\Phi(z)
\le\frac{e^{-z^2/2}}{z\sqrt{2\pi}},\qquad z>0.
\end{equation}
Consequently, by \eqref{slow:growth},
\begin{equation}\label{log:tail-size}
\log H(x)=-\beta(\log x)^2+O(\log x)
\qquad\textrm{as }x\to\infty.
\end{equation}
Moreover, for all sufficiently large \(t\), \eqref{log:mills} yields
\begin{align}\label{log:hazard}
0&<-\frac{\mathrm d}{\mathrm dt}
\log\overline\Phi\!\left(\sqrt{2\beta}(\log t-\mu)\right)\notag\\
&=\frac{\sqrt{2\beta}\exp\{-\beta(\log t-\mu)^2\}}
{\sqrt{2\pi}\,t\,\overline\Phi(\sqrt{2\beta}(\log t-\mu))}
\le C\frac{\log t}{t}.
\end{align}
Since \(t\mapsto(\log t)/t\) is eventually decreasing, it follows
that, for every fixed \(0<\theta<1\), all sufficiently large \(x\),
and every \(y\ge\theta x\),
\begin{equation}\label{log:tail-shift}
\left|\log\frac{\overline\Phi(\sqrt{2\beta}(\log y-\mu))}
{\overline\Phi(\sqrt{2\beta}(\log x-\mu))}\right|
\le C_\theta\frac{\log x}{x}|y-x|.
\end{equation}
Only the normal-tail factor is differentiated in \eqref{log:hazard};
no differentiability of \(L\) is assumed.

For every \(0<\delta<\beta\), we also have
\begin{equation}\label{log:weighted-moment}
\EE\!\left[(Z_1-1)^2
\exp\{(\beta-\delta)\log^2(1+|Z_1-1|)\}\right]<\infty.
\end{equation}
Indeed, \eqref{LogTail} and \eqref{log:tail-size} imply that
\(\PP(|Z_1-1|>t)\le\exp\{-\beta(\log t)^2+C\log t\}\)
for all sufficiently large \(t\). Integration by parts reduces the
tail part of \eqref{log:weighted-moment} to an integral bounded by a
fixed power of \(t\) times \(\exp\{-(\delta/2)(\log t)^2\}\),
which is integrable at infinity.

\begin{lemma}\label{lem:log-mgf}
Fix \(0<\delta<\beta\). There exists a constant \(C_\delta<\infty\)
such that, for \(s\in\{-1,1\}\), \(y\ge e^2\), and \(0\le\lambda\le1\)
satisfying \(\lambda y\le(\beta-\delta)(\log y)^2\),
\begin{equation}\label{log:mgf}
\EE\!\left[\exp\{\lambda s(Z_1-1)\}
\ind_{\{s(Z_1-1)\le y\}}\right]\le1+C_\delta\lambda^2.
\end{equation}
\end{lemma}

\begin{proof}
By \eqref{expz} and \eqref{Es},
\begin{align*}
&\quad\EE\!\left[\exp\{\lambda s(Z_1-1)\}
\ind_{\{s(Z_1-1)\le y\}}\right]\\
&\le1+\frac{\lambda^2}{2}\EE\!\left[(Z_1-1)^2
\exp\{\lambda[s(Z_1-1)]_+\}\ind_{\{s(Z_1-1)\le y\}}\right].
\end{align*}
For \(e^2\le u\le y\), monotonicity of \(u/(\log u)^2\) gives
\[
\frac{\lambda u}{(\log u)^2}
\le\frac{\lambda y}{(\log y)^2}\le\beta-\delta.
\]
For \(0\le u\le e^2\), we have \(e^{\lambda u}\le e^{e^2}\).
The truncated expectation above is therefore bounded by
\[
e^{e^2}\EE(Z_1-1)^2
+\EE\!\left[(Z_1-1)^2
\exp\{(\beta-\delta)\log^2(1+|Z_1-1|)\}\right].
\]
This is finite by \eqref{log:weighted-moment}, and the bound is
independent of \(s\), \(y\), and \(\lambda\). Hence
\eqref{log:mgf} follows.
\end{proof}

\begin{lemma}\label{lem:log-uniform}
Assume \eqref{LogTail}, and let \((b_n)_{n\ge1}\) satisfy
\(b_n\gg\sqrt n\log n\). Then
\begin{equation}\label{log:uniform-sum}
\lim_{n\to\infty}
\sup_{\boldsymbol s\in\{-1,1\}^n}\sup_{x\ge b_n}
\left|\frac{\PP\!\left(\sum_{k=1}^ns_k(Z_k-1)>x\right)}{nH(x)}
-\frac1n\sum_{k=1}^nc_{s_k}\right|=0.
\end{equation}
\end{lemma}

\begin{proof}
All estimates below are uniform in the deterministic signs and in
\(x\ge b_n\). The assumption on \(b_n\) implies
\(b_n/\log b_n\gg\sqrt n\), and hence
\begin{equation}\label{log:sum-boundary}
\sup_{x\ge b_n}\frac{n(\log x)^2}{x^2}\to0,\qquad
\frac{\log n}{(\log b_n)^2}\to0.
\end{equation}
Choose constants \(0<\theta<1/2\) and \(\beta<\kappa<\beta/\theta\).
We use the decomposition \eqref{Z_kPLDdec}, with this choice of
\(\theta\), and retain the notation \(R_i\) and \(E_i\) from the
proof of Lemma \ref{lem:uniform}.

\par\noindent
\textbf{(I) Control of \(P_{1n}\).}
Take \(\lambda=\kappa(\log x)^2/x\). Uniformly for \(x\ge b_n\),
\(\lambda\to0\) and
\[
\frac{\lambda\theta x}{(\log(\theta x))^2}
=\kappa\theta\frac{(\log x)^2}{(\log(\theta x))^2}
\longrightarrow\kappa\theta<\beta.
\]
Thus Lemma \ref{lem:log-mgf} applies with a fixed
\(0<\delta<\beta\) for all sufficiently large \(n\). The exponential
Markov inequality gives
\begin{equation}\label{log:P1}
P_{1n}\le\exp\{-\kappa(\log x)^2+Cn\lambda^2\}.
\end{equation}
Here \eqref{log:sum-boundary} implies
\[
\frac{n\lambda^2}{(\log x)^2}
=\kappa^2\frac{n(\log x)^2}{x^2}=o(1).
\]
Combining this estimate with \eqref{log:tail-size} and \(\kappa>\beta\),
we obtain
\begin{equation}\label{log:P1-small}
\frac{P_{1n}}{nH(x)}
\le\frac1n\exp\{-(\kappa-\beta)(\log x)^2
+Cn\lambda^2+C\log x\}=o(1).
\end{equation}
It is enough here that \(n\lambda^2=o((\log x)^2)\);
we do not require \(n\lambda^2\to0\) for this choice of \(\lambda\).

\par\noindent
\textbf{(II) Control of \(P_{2n}\).}
The identity \eqref{eq:one-jump-identity} remains valid. On
\(\{R_i\ge0\}\), the argument \(\max\{\theta x,x-R_i\}\) lies in
\([\theta x,x]\). By \eqref{slow:uniform}, the power-factor estimate
used in \eqref{Hup}, and \eqref{log:tail-shift}, we can choose a fixed
constant \(C\) such that, with \(\lambda=C(\log x)/x\),
\[
H(\max\{\theta x,x-R_i\})
\le(1+o(1))H(x)\exp\{\lambda(R_i)_+\}.
\]
Indeed, the sum of the logarithms contributed by the power and
normal-tail factors is bounded above by
\((|\gamma|/(\theta x)+C_\theta(\log x)/x)(R_i)_+\).
On \(\{R_i<0\}\), monotonicity of \(F_{s_i}\) gives
\(F_{s_i}(x-R_i)\le F_{s_i}(x)\). Consequently,
\begin{align}\label{log:one-jump}
&\quad\EE\!\left[F_{s_i}(\max\{\theta x,x-R_i\})\ind_{E_i}\right]\notag\\
&\le(c_{s_i}+\rho(\theta x))(1+o(1))H(x)
\EE\!\left[\exp\{\lambda(R_i)_+\}\ind_{E_i}\right].
\end{align}
For this choice of \(\lambda\),
\(\lambda\theta x/(\log(\theta x))^2\to0\) uniformly. Lemma
\ref{lem:log-mgf}, independence, and the argument in \eqref{Eexpup}
therefore yield
\begin{equation}\label{log:one-jump-moment}
\EE\!\left[\exp\{\lambda(R_i)_+\}\ind_{E_i}\right]
\le \exp\{Cn\lambda^2\}+\lambda\sigma\sqrt{n-1}=1+o(1),
\end{equation}
where we have used \eqref{log:sum-boundary}. Summing
\eqref{log:one-jump} over \(i\), we obtain
\begin{equation}\label{log:P2}
P_{2n}\le H(x)\left(\sum_{k=1}^n c_{s_k}+o(n)\right).
\end{equation}

\par\noindent
\textbf{(III) Control of \(P_{3n}\).}
Independence and \(F_s(t)\le CH(t)\) give
\(P_{3n}\le Cn^2H(\theta x)^2\). By \eqref{log:tail-size},
\begin{equation}\label{log:P3}
\frac{P_{3n}}{nH(x)}
\le Cn\frac{H(\theta x)^2}{H(x)}
\le C\exp\{\log n-\beta(\log x)^2+C\log x\}=o(1),
\end{equation}
where the last assertion follows from \eqref{log:sum-boundary}.
Combining \eqref{log:P1-small}, \eqref{log:P2}, and \eqref{log:P3},
we obtain
\begin{equation}\label{log:sum-upper}
\PP\!\left(\sum_{k=1}^ns_k(Z_k-1)>x\right)
\le H(x)\left(\sum_{k=1}^n c_{s_k}+o(n)\right).
\end{equation}

Next, we give the lower bound. For all sufficiently large \(n\), take
\(d_n=n^{1/4}(b_n/\log b_n)^{1/2}\). Then
\[
\sqrt n\ll d_n\ll b_n/\log b_n.
\]
The disjoint events used in \eqref{PLDlow} are still contained in
\(\{\sum_{k=1}^ns_k(Z_k-1)>x\}\), and \eqref{PR_ilow} gives
\begin{equation}\label{log:lower-event}
\PP(R_i\ge-d_n,E_i)
\ge1-CnH(\theta x)-\frac{n\sigma^2}{d_n^2}=1-o(1).
\end{equation}
Here \(nH(\theta x)\to0\) uniformly by \eqref{log:tail-size} and
\eqref{log:sum-boundary}. Since \(d_n/x\to0\) uniformly,
\eqref{slow:uniform} gives
\[
\frac{(x+d_n)^\gamma L(x+d_n)}{x^\gamma L(x)}\to1.
\]
The normal-tail factor satisfies, by \eqref{log:tail-shift},
\[
\left|\log\frac{\overline\Phi(\sqrt{2\beta}(\log(x+d_n)-\mu))}
{\overline\Phi(\sqrt{2\beta}(\log x-\mu))}\right|
\le C\frac{d_n\log x}{x}
\le C\frac{d_n\log b_n}{b_n}\to0.
\]
Thus \(H(x+d_n)/H(x)\to1\) uniformly for \(x\ge b_n\), and
\[
\sum_{i=1}^n F_{s_i}(x+d_n)
=H(x)\left(\sum_{i=1}^n c_{s_i}+o(n)\right).
\]
Using independence of \(s_i(Z_i-1)\) and \((R_i,E_i)\), together
with \eqref{log:lower-event}, we conclude that
\begin{align}\label{log:sum-lower}
\PP\!\left(\sum_{k=1}^ns_k(Z_k-1)>x\right)
&\ge\sum_{i=1}^n F_{s_i}(x+d_n)\PP(R_i\ge-d_n,E_i)\notag\\
&\ge H(x)\left(\sum_{k=1}^n c_{s_k}-o(n)\right).
\end{align}
Combining \eqref{log:sum-upper} and \eqref{log:sum-lower} completes
the proof of Lemma \ref{lem:log-uniform}.
\end{proof}
\section{Proof of Main Results}
\subsection{Proof of Theorem \ref{thm:subexponential}}
We first fix \(n\). Since \(|T_n|\le n\) and
\(S_n=T_n+\sum_{k=1}^nX_k(Z_k-1)\), subexponentiality and the
long-tail property give
\[
\PP(|S_n|>x)
\le\PP\!\left(\sum_{k=1}^n|Z_k-1|>x-n\right)
\sim n\overline F(x)
\qquad\textrm{as }x\to\infty.
\]
For the lower bound, the events
\[
\left\{|Z_i-1|>x+n+\sum_{j\ne i}|Z_j-1|\right\},
\qquad i=1,\ldots,n,
\]
are pairwise disjoint. On each of these events, the reverse triangle
inequality implies that
\[
|S_n|\ge |Z_i-1|-\sum_{j\ne i}|Z_j-1|-n>x.
\]
Consequently, by independence of the random step sizes,
\begin{align*}
\PP(|S_n|>x)
&\ge\sum_{i=1}^n
\PP\!\left(|Z_i-1|>x+n+\sum_{j\ne i}|Z_j-1|\right)\\
&=\sum_{i=1}^n
\EE\!\left[\overline F\!\left(x+n+\sum_{j\ne i}|Z_j-1|\right)\right].
\end{align*}
For fixed \(n\) and \(i\), the long-tail property yields
\[
\frac{\overline F(x+n+\sum_{j\ne i}|Z_j-1|)}{\overline F(x)}
\to1\qquad\textrm{a.s}.
\]
The ratio is between \(0\) and \(1\). Thus, by the dominated
convergence theorem, the preceding lower bound is asymptotic to
\(n\overline F(x)\). Combining the upper and lower bounds, we obtain
\[
\PP(|S_n|>x)\sim n\overline F(x)
\qquad\textrm{as }x\to\infty
\quad\textrm{for every fixed }n.
\]

When \(q=1/2\), the conditional distribution of \(X_k\) in the
Introduction gives
\((X_1,\ldots,X_n)\overset{d}=(-X_1,\ldots,-X_n)\).
Since the step-size sequence is independent of the direction sequence,
\(S_n\overset{d}=-S_n\). Therefore,
\begin{equation}\label{subexp:fixed-tail}
\PP(S_n>x)=\PP(S_n<-x)=\frac12\PP(|S_n|>x)
\sim\frac n2\overline F(x)
\end{equation}
as \(x\to\infty\), for every fixed \(n\).
For each \(n\), choose \(b_n\ge n\) sufficiently large such that
\[
\sup_{x\ge b_n}
\left|\frac{\PP(S_n>x)}{n\overline F(x)}-\frac12\right|
\le\frac1n.
\]
Then \(b_n\to\infty\), and \eqref{subexp:main-right} follows.
The left-tail conclusion follows by symmetry.

Finally, for \(x>1\),
\[
\overline F(x+1)\le\PP(|Z_1|>x)\le\overline F(x-1),
\]
so the long-tail property gives
\(\PP(|Z_1|>x)\sim\overline F(x)\).
The inequalities
\[
\PP(S_n>x)\le\PP(S_n\ge x)\le\PP(S_n>x-1)
\]
and the analogous inequalities for the left tail and the step-size
tails show that all the corresponding fixed-\(n\) limits in
\eqref{subexp:fixed-tail} remain valid. Increasing \(b_n\), if
necessary, so that the finitely many corresponding ratios satisfy
the same bound \(1/n\), proves the remaining assertions.
Hence, we complete the proof of Theorem \ref{thm:subexponential}.
\subsection{Proof of Theorem \ref{thm:main}}\label{sec:weibull-proof}
\quad We retain the assumptions and notation of Theorem \ref{thm:main} and the
Preliminary lemmas, with \(H\) as in \eqref{H:weibull}, $\EE Z_1=1$, and $q=1/2$.
We first give the proof of the right-tail conclusion. 

For $0<p\le3/4$, let
\[
a_k=\frac{\Gamma(k)\Gamma(2p)}{\Gamma(k+2p-1)},\qquad k\ge1.
\]
For $n\ge3$, define
\begin{equation}\label{proof:variance-scale}
v_n=
\begin{cases}
\displaystyle\frac{\sum_{k=3}^n(k-1)^2}{(n-1)^2},&p=0,\\[6pt]
\displaystyle\frac{1}{a_n^2}\sum_{k=1}^n a_k^2,&0<p\le3/4.
\end{cases}
\end{equation}
Here $v_n$ is the deterministic scale for $T_n$ itself. For $p>0$,
the gamma-ratio asymptotic gives $a_n\sim\Gamma(2p)n^{1-2p}$.
Summing $a_k^2$, and using direct summation when $p=0$, we have
\begin{equation}\label{proof:variance-asymptotics}
v_n\sim
\begin{cases}
\displaystyle\frac{n}{3-4p},&0\le p<3/4,\\[5pt]
n\log n,&p=3/4.
\end{cases}
\end{equation}
Thus, by \eqref{proof:boundary}, it holds that
\begin{equation}\label{proof:small}
v_n b_n^{-2(1-\alpha)}\to0,\qquad
\sqrt{v_n}\,b_n^{\alpha-1}\to0,\qquad
v_n^{-1}b_n^{2-\alpha}\to\infty.
\end{equation}
Moreover, $\EE T_n=0$, and the identity
$\sum_{k=1}^n\ind_{\{X_k=1\}}=(n+T_n)/2$ gives
\[
\sum_{k=1}^n c_{X_k}
=\frac{n+T_n}{2}c_1+\frac{n-T_n}{2}c_{-1}.
\]
Consequently,
\begin{equation}\label{proof:coefficient}
0\le\frac1n\sum_{k=1}^n c_{X_k}
\le\max\{c_1,c_{-1}\},\qquad
\frac1n\EE\!\left[\sum_{k=1}^n c_{X_k}\right]
=\frac{c_1+c_{-1}}2.
\end{equation}

\par
An Azuma--Hoeffding argument, as used in the martingale approach to
ERW deviation inequalities (see Hu and Fan \cite{HFDI}), yields
\begin{equation}\label{proof:one-sided}
\PP(T_n>u)=\PP(T_n<-u)
\le\exp\!\left\{-\frac{u^2}{2v_n}\right\},\qquad u\ge0.
\end{equation}
The equality of the two tails follows from the symmetry of $T_n$ when $q=1/2$.
For completeness, when $0<p\le3/4$, $(a_kT_k)_{k\ge1}$ is a martingale and
\[
a_kT_k-a_{k-1}T_{k-1}
=a_k\left(X_k-\frac{2p-1}{k-1}T_{k-1}\right),\qquad k\ge2.
\]
Conditional on $\GG_{k-1}$, this increment has mean zero and takes values
in an interval of length $2a_k$. Also, $a_1T_1=T_1$ has mean zero and
takes values in an interval of length $2$. By the conditional Hoeffding
inequality and iterated expectations,
\[
\EE e^{\lambda T_n}\le\exp\{\lambda^2v_n/2\},\qquad\lambda\in\mathbf R.
\]
When $p=0$, we have $T_2=0$, and $((k-1)T_k)_{k\ge2}$ is a martingale
whose increment at time $k\ge3$ is $(k-1)X_k+T_{k-1}$. Its conditional
mean is zero and its conditional range has length $2(k-1)$, so the same
exponential-moment bound holds with the $p=0$ scale in
\eqref{proof:variance-scale}. The exponential Markov inequality, with
$\lambda=u/v_n$, proves \eqref{proof:one-sided}.
Consequently,
\begin{equation}\label{proof:concentration}
\PP(|T_n|>u)
\le2\exp\!\left\{-\frac{u^2}{2v_n}\right\},\qquad u\ge0.
\end{equation}
Applying the one-sided bound \eqref{proof:one-sided} at $u=x/2$ and
using \eqref{slow:growth} and \eqref{proof:small}, we can deduce that
\begin{equation}\label{proof:exception}
\sup_{x\ge b_n}
\frac{\PP(T_n>x/2)}{nH(x)}
\le\sup_{x\ge b_n}\frac{1}{nx^\gamma L(x)}
\exp\!\left\{\beta x^\alpha-\frac{x^2}{8v_n}\right\}\to0.
\end{equation}
Indeed, $x^{2-\alpha}/v_n\ge b_n^{2-\alpha}/v_n\to\infty$, so the
negative term in the exponent dominates uniformly for $x\ge b_n$.
The factor \(1/(x^\gamma L(x))\) is bounded by a power of \(x\) by
\eqref{slow:growth} and does not affect this limit.

\par

Conditional on $\GG_n$, the variables $(X_k(Z_k-1))_{1\le k\le n}$ are
independent centered variables of the form treated in Lemma \ref{lem:uniform}.
The sequence \(b_n/2\) satisfies the lower-boundary condition in that lemma.
In particular, the random threshold
$\max\{x/2,x-T_n\}$ is $\GG_n$-measurable and is at least $b_n/2$.
Therefore, it holds that
\begin{align}\label{proof:conditional}
&\quad\PP\!\left(S_n-T_n>\max\{x/2,x-T_n\}\mid\GG_n\right)\notag\\
&=H\!\left(\max\{x/2,x-T_n\}\right)
\left(\sum_{k=1}^n c_{X_k}+o(n)\right).
\end{align}
The error \(o(n)\) in \eqref{proof:conditional}, divided by \(n\),
is bounded uniformly over all realizations of the signs and all
\(x\ge b_n\) by a deterministic sequence tending to zero.

We next estimate the effect of the random shift on both factors of \(H\).
By the mean value theorem, we have
\begin{equation}\label{proof:shift-power}
\left|\max\{x/2,x-r\}^{\alpha}-x^\alpha\right|
\le\alpha(x/2)^{\alpha-1}|r|,
\end{equation}
where $r\in\mathbf R$.
Similarly, since \(\max\{x/2,x-r\}\ge x/2\) and
\(\left|\max\{x/2,x-r\}-x\right|\le|r|\), it holds that
\begin{equation}\label{proof:shift-log}
\left|\log\frac{\max\{x/2,x-r\}}{x}\right|
\le\frac{2|r|}{x}.
\end{equation}
Fix \(\varepsilon>0\). By \eqref{rv:potter} and
\eqref{proof:shift-log}, for all sufficiently large \(n\), uniformly
in \(x\ge b_n\) and \(r\in\mathbf R\),
\begin{equation}\label{proof:shift-L}
\left|\log\frac{\max\{x/2,x-r\}^{\gamma}L(\max\{x/2,x-r\})}
{x^\gamma L(x)}\right|
\le\log(1+\varepsilon)+\frac{2(|\gamma|+1)|r|}{x}.
\end{equation}
Take \(\lambda=\beta\alpha(x/2)^{\alpha-1}+2(|\gamma|+1)/x\). By \eqref{proof:small},
\begin{equation}\label{proof:lambda}
\sup_{x\ge b_n}v_n\lambda^2\to0.
\end{equation}
Using $|e^z-1|\le e^{|z|}-1$, \eqref{proof:shift-power}, and
\eqref{proof:shift-L}, we obtain
\begin{equation}\label{proof:shift-bound}
\left|
\frac{H(\max\{x/2,x-T_n\})}{H(x)}-1
\right|
\le (1+\varepsilon)\exp\{\lambda|T_n|\}-1.
\end{equation}
Moreover, by Tonelli's theorem, \eqref{proof:concentration} and \eqref{proof:lambda}, it holds that
\begin{align}\label{proof:shift-moment}
\EE(\exp\{\lambda|T_n|\}-1)
&=\lambda\int_0^\infty e^{\lambda u}\PP(|T_n|>u)\,du
\le2\lambda\int_0^\infty
       \exp\!\left\{\lambda u-\frac{u^2}{2v_n}\right\}\,du\notag\\
&=2\lambda\sqrt{v_n}\,e^{v_n\lambda^2/2}
       \int_{-\lambda\sqrt{v_n}}^\infty e^{-t^2/2}\,dt
\le2\sqrt{2\pi}\,\lambda\sqrt{v_n}\,e^{v_n\lambda^2/2}
=o(1).
\end{align}
Consequently, the limit superior of the expectation on the left-hand
side of \eqref{proof:shift-bound}, uniformly for \(x\ge b_n\), is at
most \(\varepsilon\). Letting \(\varepsilon\downarrow0\), we have
\begin{equation}\label{proof:shift-average}
\sup_{x\ge b_n}\EE\!\left[
\left|\frac{H(\max\{x/2,x-T_n\})}{H(x)}-1\right|
\right]\to0.
\end{equation}
By \eqref{proof:coefficient} and \eqref{proof:shift-average}, the error
after multiplication by $\sum_{k=1}^n c_{X_k}$ is $o(n)$ in expectation.
It follows that
\begin{align}\label{proof:averaged-tail}
&\quad\EE\!\left[
\left(\sum_{k=1}^n c_{X_k}\right)
H(\max\{x/2,x-T_n\})\right]\notag\\
&=H(x)\left\{\EE\!\left[\sum_{k=1}^n c_{X_k}\right]+o(n)\right\}\notag\\
&=\frac{n(c_1+c_{-1})H(x)}2(1+o(1)).
\end{align}
By \eqref{proof:shift-average}, we also have
\(\EE H(\max\{x/2,x-T_n\})=H(x)(1+o(1))\) uniformly for \(x\ge b_n\).
Thus the additive error in \eqref{proof:conditional} contributes
\(o(nH(x))\) after taking expectations. Using
\eqref{proof:averaged-tail}, we obtain
\begin{equation}\label{proof:restricted-tail}
\PP\!\left(S_n-T_n>\max\{x/2,x-T_n\}\right)
=\frac{n(c_1+c_{-1})H(x)}2(1+o(1)).
\end{equation}

The event on the left-hand side of \eqref{proof:restricted-tail} is
exactly $\{S_n>x,\ S_n-T_n>x/2\}$. Since
\[
\{S_n>x,\ S_n-T_n\le x/2\}\subseteq\{T_n>x/2\},
\]
we have
\begin{align}\label{proof:comparison}
0&\le\PP(S_n>x)
-\PP\!\left(S_n-T_n>\max\{x/2,x-T_n\}\right)\notag\\
&\le\PP(T_n>x/2)=o(nH(x)).
\end{align}
Combining \eqref{proof:restricted-tail}, \eqref{proof:comparison}, and
$\overline F(x)\sim(c_1+c_{-1})H(x)$, we can deduce that
\begin{equation}\label{proof:right-tail}
\PP(S_n>x)
=\frac{n(c_1+c_{-1})H(x)}2(1+o(1))
=\frac n2\overline F(x)(1+o(1))
\end{equation}
uniformly for $x\ge b_n$. Hence, the right-tail conclusion follows.
When \(q=1/2\), the joint distribution of \((X_1,\ldots,X_n)\) is
invariant under simultaneous sign reversal. Independence of the step
sizes therefore implies that \(S_n\) and \(-S_n\) have the same
distribution, which proves the left-tail conclusion.

\par\noindent

For $x>1$, the triangle inequality gives
\begin{equation}\label{proof:step-shift}
\overline F(x+1)\le\PP(|Z_1|>x)\le\overline F(x-1).
\end{equation}
Since $\alpha<1$ and $b_n\to\infty$, the mean value theorem and
\eqref{slow:uniform} yield
\[
\sup_{x\ge b_n}\left|\frac{H(x\pm1)}{H(x)}-1\right|\to0.
\]
Together with $\overline F(t)\sim(c_1+c_{-1})H(t)$, this implies that
\[
\sup_{x\ge b_n}
\left|\frac{\PP(|Z_1|>x)}{\overline F(x)}-1\right|\to0.
\]
Consequently, $\overline F(x)$ can be replaced by $\PP(|Z_1|>x)$ in
both tail conclusions.

Finally, we have
\[
\PP(S_n>x)\le\PP(S_n\ge x)\le\PP(S_n>x-1).
\]
The lower boundary $b_n/2$ also satisfies \eqref{proof:boundary}, so
\eqref{proof:right-tail} applies at $x-1$ for every $x\ge b_n$ and all
sufficiently large $n$. Since $\overline F(x-1)/\overline F(x)\to1$
uniformly on this range, the same precise large deviation holds for
$\PP(S_n\ge x)$. By symmetry, it also holds for $\PP(S_n\le-x)$.
If non-strict step-size tails are used in the denominator, the
inequalities
\[
\overline F(x)\le\PP(|Z_1-1|\ge x)\le\overline F(x-1)
\]
and the corresponding inequalities for $|Z_1|$ give the same conclusion.
This completes the proof of Theorem \ref{thm:main} under \eqref{proof:boundary}.
\qed

\subsection{Proof of Theorem \ref{thm:lognormal}}\label{sec:lognormal-proof}
\quad We now retain the assumptions of Theorem \ref{thm:lognormal},
with \(H\) as in \eqref{H:lognormal}. The scale \(v_n\), the coefficient
identity \eqref{proof:coefficient}, and the concentration inequalities
\eqref{proof:one-sided} and \eqref{proof:concentration} are unchanged.
By \eqref{proof:variance-asymptotics} and \eqref{log:boundary},
\begin{equation}\label{log:erw-small}
\frac{v_n(\log b_n)^2}{b_n^2}\to0,\qquad
\sup_{x\ge b_n}\frac{v_n(\log x)^2}{x^2}\to0.
\end{equation}
In particular, \(b_n/2\gg\sqrt n\log n\), as required in
Lemma \ref{lem:log-uniform}.

We first control the exceptional event \(\{T_n>x/2\}\). By
\eqref{proof:one-sided} and \eqref{log:tail-size},
\begin{equation}\label{log:exception}
\sup_{x\ge b_n}\frac{\PP(T_n>x/2)}{nH(x)}
\le\sup_{x\ge b_n}\frac1n
\exp\!\left\{\beta(\log x)^2+C\log x-\frac{x^2}{8v_n}\right\}
\to0.
\end{equation}
Indeed, \eqref{log:erw-small} implies that
\(x^2/(v_n(\log x)^2)\to\infty\) uniformly for \(x\ge b_n\), so
the negative term dominates uniformly.

Conditional on \(\GG_n\), the variables \((X_k(Z_k-1))_{1\le k\le n}\)
are independent and centered. Applying Lemma \ref{lem:log-uniform}
with lower boundary \(b_n/2\) gives
\begin{align}\label{log:conditional}
&\quad\PP\!\left(S_n-T_n>\max\{x/2,x-T_n\}\mid\GG_n\right)\notag\\
&=H(\max\{x/2,x-T_n\})
\left(\sum_{k=1}^n c_{X_k}+o(n)\right),
\end{align}
where the error divided by \(n\) is bounded by a deterministic
sequence tending to zero, uniformly over all sign realizations and
all \(x\ge b_n\).

Next, we estimate the effect of the random shift. Fix
\(\varepsilon>0\). By \eqref{rv:potter}, \eqref{log:tail-shift} with
\(\theta=1/2\), and \eqref{proof:shift-log}, we have, for all
sufficiently large \(n\),
\begin{align}\label{log:shift-bound}
&\quad\left|\log\frac{H(\max\{x/2,x-T_n\})}{H(x)}\right|\notag\\
&\le\log(1+\varepsilon)
+\left(C\frac{\log x}{x}+\frac{2(|\gamma|+1)}x\right)|T_n|,
\qquad x\ge b_n.
\end{align}
Take \(\lambda=C(\log x)/x+2(|\gamma|+1)/x\). Then
\eqref{log:erw-small} yields
\[
\sup_{x\ge b_n}v_n\lambda^2\to0.
\]
The integration argument in \eqref{proof:shift-moment} therefore
implies that
\begin{equation}\label{log:shift-moment}
\sup_{x\ge b_n}\EE(e^{\lambda|T_n|}-1)
\le\sup_{x\ge b_n}
2\sqrt{2\pi}\,\lambda\sqrt{v_n}\,e^{v_n\lambda^2/2}\to0.
\end{equation}
Using \(|e^z-1|\le e^{|z|}-1\) in \eqref{log:shift-bound}, we obtain
\[
\EE\!\left|\frac{H(\max\{x/2,x-T_n\})}{H(x)}-1\right|
\le\varepsilon+(1+\varepsilon)\EE(e^{\lambda|T_n|}-1).
\]
Taking the supremum over \(x\ge b_n\), then letting \(n\to\infty\)
and \(\varepsilon\downarrow0\), gives
\begin{equation}\label{log:shift-average}
\sup_{x\ge b_n}\EE\!\left|
\frac{H(\max\{x/2,x-T_n\})}{H(x)}-1\right|\to0.
\end{equation}

By \eqref{proof:coefficient} and \eqref{log:shift-average},
\begin{align}\label{log:averaged-tail}
&\quad\EE\!\left[
H(\max\{x/2,x-T_n\})\sum_{k=1}^n c_{X_k}\right]\notag\\
&=H(x)\left\{\EE\!\left[\sum_{k=1}^n c_{X_k}\right]+o(n)\right\}
=\frac{n(c_1+c_{-1})}{2}H(x)(1+o(1)).
\end{align}
Also, \eqref{log:shift-average} implies that
\(\EE H(\max\{x/2,x-T_n\})=H(x)(1+o(1))\), so the additive error
in \eqref{log:conditional} contributes \(o(nH(x))\) after taking
expectations. Consequently,
\begin{equation}\label{log:restricted-tail}
\PP\!\left(S_n-T_n>\max\{x/2,x-T_n\}\right)
=\frac{n(c_1+c_{-1})}{2}H(x)(1+o(1)).
\end{equation}
As in \eqref{proof:comparison}, the difference between this probability
and \(\PP(S_n>x)\) satisfies
\begin{align}\label{log:comparison}
0&\le\PP(S_n>x)
-\PP\!\left(S_n-T_n>\max\{x/2,x-T_n\}\right)\notag\\
&\le\PP(T_n>x/2)=o(nH(x)).
\end{align}
Combining \eqref{log:restricted-tail}, \eqref{log:comparison}, and
\(\overline F(x)\sim(c_1+c_{-1})H(x)\), we obtain
\begin{equation}\label{log:right-tail}
\PP(S_n>x)
=\frac{n(c_1+c_{-1})}{2}H(x)(1+o(1))
=\frac n2\overline F(x)(1+o(1))
\end{equation}
uniformly for \(x\ge b_n\). The joint sign symmetry when \(q=1/2\),
as used in Section \ref{sec:weibull-proof}, gives the left-tail
conclusion.

Finally, \eqref{slow:uniform} and \eqref{log:tail-shift} imply
\[
\sup_{x\ge b_n}\left|\frac{H(x\pm1)}{H(x)}-1\right|\to0.
\]
The comparison \eqref{proof:step-shift} therefore allows
\(\overline F(x)\) to be replaced by \(\PP(|Z_1|>x)\).
Moreover,
\[
\PP(S_n>x)\le\PP(S_n\ge x)\le\PP(S_n>x-1).
\]
The sequence \(b_n/2\) also satisfies \eqref{log:boundary}, so
\eqref{log:right-tail} applies at \(x-1\) for every \(x\ge b_n\)
and all sufficiently large \(n\). Since
\(\overline F(x-1)/\overline F(x)\to1\) uniformly, the non-strict
right-tail conclusion follows. The non-strict left-tail conclusion
follows by symmetry. The inequalities
\[
\overline F(x)\le\PP(|Z_1-1|\ge x)\le\overline F(x-1)
\]
and their counterparts for \(|Z_1|\) give the corresponding
conclusions for non-strict step-size tails. This completes the proof
of Theorem \ref{thm:lognormal}.
\qed


\begin{thebibliography}{99}
\bibitem{BMERW} Bercu B. A martingale approach for the elephant random walk. J Phys A: Math Theor, 2018, 51: 015201.
\bibitem{B21} Bertoin J. Counting the zeros of an elephant random walk. Trans Amer Math Soc, 2022, 375: 5539--5560.
\bibitem{BGT} Bingham N H, Goldie C M, Teugels J L. Regular Variation. Cambridge: Cambridge University Press, 1987.
\bibitem{CGSCLT} Coletti C F, Gava R, Schütz G M. Central limit theorem and related results for the elephant random walk. J Math Phys, 2017, 58: 053303.
%\bibitem{DFHRCCLT} Dedecker J, Fan X, Hu H, et al. Rates of convergence in the central limit theorem for the elephant random walk with random step sizes. J Stat Phys, 2023, 190: 154.
\bibitem{DDS}Denisov D, Dieker A B, Shneer V. Large deviations for random walks under subexponentiality: the big-jump domain. Ann Probab, 2008, 36: 1946--1991.
\bibitem{FanHuMa2021} Fan X, Hu H, Ma X. Cram\'er moderate deviations for the elephant random walk. J Stat Mech Theory Exp, 2021, 2: 023402.
\bibitem{FanLiangMa} Fan X, Liang S, Ma D. Cram\'er-type moderate deviations and precise large deviations for the elephant random walk with random step sizes. Sci China Math, Advance online publication. DOI: 10.1007/s11425-026-2671-6.
\bibitem{FSCMD} Fan X, Shao Q M. Cram\'er's moderate deviations for martingales with applications. Ann Inst H Poincaré Probab Statist, 2024, 60: 2046--2074.
\bibitem{HFDI} Hu H, Fan X. Deviation inequalities for the elephant random walk with random step sizes. J Math Phys, 2025, 66: 033301.
\bibitem{Liang2026} Liang S. Local precise large deviations for the elephant random walk with random step sizes. Stat Probab Lett, 2026, 239: 110905.
\bibitem{NAVILTLDP} Nagaev A V. Integral limit theorems for large deviations when Cram\'er's condition is not fulfilled I, II. Theory Probab Appl, 1969, 14: 51--64, 193--208.
\bibitem{NSVLDP} Nagaev S V. Large deviations of sums of independent random variables. Ann Probab, 1979, 7: 745--789.
\bibitem{schutz2004elephants} Schütz G M, Trimper S. Elephants can always remember: exact long-range memory effects in a non-markovian random walk. Phys Rev E, 2004, 70: 045101.
\end{thebibliography}
\end{document}